# High Order Riesz Transforms and Mean Value Formula for Generalized Translate Operator

I. Ekincioglu, H.H. Sayan and C. Keskin

**Abstract.** In this paper, the mean value formula depends on the Bessel generalized shift operator corresponding to the solutions of the boundary value problem related to the Bessel operator are studied. In addition to, Riesz Bessel transforms $R_{\mathrm{B}_{x_j}}$ related to the Bessel operators are studied. Since Bessel generalized shift operator is translation operator corresponding to the Bessel operator, we construct a family of $R_{\mathrm{B}_{x_j}}$ by using Bessel generalized shift operator. Finally, we analysis weighted inequalities involving $R_{\mathrm{B}_{x_j}}$.

## 1. Introduction

Singular integral operators are playing an important role in Harmonic Analysis, theory of functions and partial differential equations. Singular integrals associated with the $\Delta_{\mathrm{B}_n}$ Laplace-Bessel differential operator which is known as an important operator in analysis and its applications, have been the research areas of many mathematicians such as B. Muckenhoupt and E.M. Stein [17, 19], I. Kipriyanov and M. Klyuchantsev [10],[12], I.A. Aliev and A.D. Gadjiev [2], V.S. Guliyev [7], A.D. Gadjiev . Also, singular integral operators related to generalized shift operator were studied by [1] and [4, 5, 6] and others.

The Bessel generalized shift operator is one of the most important generalized shift operator on the half line $\mathbb{R}_+$, [14, 8]. The Bessel generalized translation is used while studying various problems connected with Bessel operators [9]. Fourier-Bessel harmonic analysis, i.e. the part of harmonic analysis addressing various problems on Bessel (Hankel) integral transforms, is closely connected with Bessel generalized shift operator.

It is well known that the fundamental solutions of the classical stationary of mathematical physics ( the harmonic equation, polyharmonic equation and Helmholtz equation) are radial functions. Therefore, it is natural to seek these solutions as solutions of an ordinary differential equations. However, since the spherical coordinate transformation brings an equation with the Laplace operator in $\mathbb{R}^n$ to an ordinary differential equation with the singular Bessel differential operator, an interest arose (probably, a long time ago) in studying methods for constructing the fundamental solutions of singular ordinary differential equations with the Bessel operator in place of the second derivative. In this connection, it might be very useful to prove a theorem on the fundamental solution of an ordinary differential equation involving the Bessel operator with constant coefficients, similar to the well-known theorem on the fundamental solution of an ordinary differential equation.

The obtained result on the fundamental solution of an ordinary differential operator with the Bessel operator has allowed us to analyze equations with the singular differential operator $\triangle_{\mathrm{B}_\gamma} = \sum_{i=1}^{n} \mathrm{B}_{\gamma_i} + \sum_{i=n+1}^{N} \dfrac{\partial^2}{\partial x_i^2}$ where the different indices $\gamma_i$ act with respect to part of the variables and these indices may take the negative values. The latter fact is essential because the tools used in problems of this kind (The Poisson operator and generalized shift operator of integral nature) are defined only for $\gamma_i > 0$. To find the fundamental solution with a singularity at an arbitrary point, we use generalized shift operator that acts with respect to the radial variable. Note also that the mixed-type generalized shift operator, which is conventionally used in such problems, coincides with the radial shift operator on radial functions provided that $\gamma_i > 0$. The found fundamental solutions(of the B-harmonic and B-polyharmonic equations and of the singular Helmholtz equation) coincide with the known solutions when $\gamma_i > 0$ and with the classical solutions when $\gamma_i = 0$.

It is well known that harmonic functions satisfy various mean value theorems which may be considered as generalizations of the Gauss mean value theorem. There have been a number of studies on mean value theorems. Cheng [3] obtained a converse for a different mean value expression. Nicolesco [18] gave an expression in terms of certain iterated means and showed that a converse was also true. The mean value theorems for harmonic functions has also been carried out by Pizetti [23], Picone [22], Ekincioglu [4] and Kipriyanova [11] and [13].

The solutions of the boundary value problems for Laplace operator are related to the ordinary shift operator. Also, the solutions of boundary value problems for Laplace-Bessel and Bessel operator are corresponding to the generalized shift operator and Bessel generalized shift operator, respectively.

In this paper, singular integral operators generated by Bessel generalized shift operator are studied. In addition, the mean value formulas related to Bessel generalized shift operator are given.

Riesz Bessel singular integral operators related to generalized shift operator for Laplace-Bessel operator were showed in [1] and [5]. The authors used the mean value theorem related to the generalized shift operator. The mean value formula for equations $\Delta_{\mathrm{B}}u = 0$ and $\Delta_{\mathrm{B}}u + \lambda u = 0$ were obtained by [4] and [13], respectively.

In this study, we introduce the mean value formula for equation $\mathrm{B}_{x_i}u = 0$ and the high order Riesz Bessel transform associated with Bessel generalized shift operator for $\mathrm{B}_{x_i}$ Bessel differential operator

$$\mathrm{B}_{x_i}f = \sum_{i=1}^{n} \partial_{x_i}^2 + \frac{2\gamma_j}{x_i}\partial_{x_i} \quad , \gamma_i > 0$$

where $\partial_{x_i}^2 = \dfrac{\partial^2}{\partial x_i^2}$ .

Let $\mathbb{R}^n$ be $n$-dimensional Euclidean space and $x = (x_1, ..., x_n)$, $\xi = (\xi_1, ..., \xi_n)$ be vectors in $\mathbb{R}^n$, then $x \cdot \xi = x_1\xi_1 + \ldots + x_n\xi_n$, $|x| = (x \cdot x)^{1/2}$. Denote $\mathbb{R}_+^n = \{x \in \mathbb{R}^n; x_1 > 0, \ldots, x_k > 0,\ 1 \le k \le n\}$, $S_+^n = \{x \in \mathbb{R}_+^n : |x| = 1\}$, $\gamma = (\gamma_1, \ldots, \gamma_n)$, $\gamma_1 > 0, \ldots, \gamma_n > 0$, $|\gamma| = \gamma_1 + \ldots + \gamma_n$ and $d\mu_\gamma(x) = \prod_{i=1}^{n} x_i^{2\gamma_i} dx$. We shall denote by $L_{p,\gamma}(\mathbb{R}^n, d\mu_\gamma(x))$-spaces (The Lebesgue space with respect to the measure $\mu_\gamma$), the set of all measurable functions $f$ on $\mathbb{R}_+^n$ such that the norm

$$\|f\|_{L_{p,\gamma}} \equiv \|f\|_{p,\gamma} = \left( \int_{\mathbb{R}_+^n} |f(x)|^p d\mu_\gamma(x) \right)^{1/p}, \quad 1 \le p < \infty$$

is finite. We determine the Bessel generalized translation $T^y\varphi(x) = u(x)$, $x, y \in \mathbb{R}_+^n$ of a function $\varphi(x) \in \mathcal{C}^{(2)}(\mathbb{R}_+^n)$ as the solution to the following initial value problem:

$$\left.\begin{array}{l} \mathrm{B}_{x_i}u(x,y) = \mathrm{B}_{y_i}u(x,y), \\ u(x,0) = \varphi(x), \quad u_y(x,0) = 0, \end{array}\right\} \tag{1.1}$$

where $i = 1, 2, \ldots, n$ and $\mathrm{B}_{x_i}$ is the Bessel differential operator. The solution of the initial value problem (1.1) exists, is unique, and can be written explicitly as

$$\begin{aligned} T^y\varphi(x) \quad &:= c_\gamma \int_0^\pi \ldots \int_0^\pi \varphi\big((x_1, y_1)_{\alpha_1}, \ldots, (x_n, y_n)_{\alpha_n}\big) \\ &\times \Big( \prod_{i=1}^{n} \sin^{2\gamma_i - 1} \alpha_i \Big) d\alpha_1 \ldots d\alpha_n, \end{aligned} \tag{1.2}$$

where $c_\gamma = \prod_{i=1}^{n} \Gamma(\gamma_i + \frac{1}{2}) \big[\Gamma(\frac{1}{2})\Gamma(\gamma_i)\big]^{-1}$ and $(x_i, y_i)_{\alpha_i} = \sqrt{x_i^2 + y_i^2 - 2x_i y_i \cos\alpha_i}$, $1 \le i \le n$. By formula (1.2), the operator $T^y$ can be extended to all functions $L_{p,\gamma}(\mathbb{R}_+^n)$. The operator $T^y$ satisfying (1.1) may be regarded as a Bessel-generalized shift operator Bessel Generalized Shift Operator (see [12], [14] and [15]). We remark that this shift operator is closely connected with the

Bessel differential operator. The convolution operator determined by $T^y$ is as follows:

$$(f * \varphi)(x) = \int_{\mathbb{R}^n_+} f(y) T^y \varphi(x) d\mu_\gamma(y). \tag{1.3}$$

Convolution (1.3) is known as a B-convolution. We note some properties for the B-convolution and the Bessel generalized shift operator:

- If $f(x), \varphi(x) \in \mathcal{C}(\mathbb{R}^n_+)$, $\varphi(x)$ is a bounded function, $x > 0$ and

$$\int_0^\infty |f(x)| d\mu_\gamma(x) < \infty$$

  then

$$\int_{\mathbb{R}^n_+} T^y f(x) \varphi(y) d\mu_\gamma(y) = \int_{\mathbb{R}^n_+} f(y) T^y \varphi(x) d\mu_\gamma(y).$$

- From the above result, we have the following equality for $\varphi(x) = 1$.

$$\int_{\mathbb{R}^n_+} T^y f(x) d\mu_\gamma(y) = \int_{\mathbb{R}^n_+} f(y) d\mu_\gamma(y).$$

- $(f * \varphi)(x) = (\varphi * f)(x)$.

The Fourier-Bessel transform is defined and invertible on functions $\varphi \in \mathcal{S}(\mathbb{R}^n_+)$,

$$[F_{\mathrm{B}}\varphi](y) = c_\gamma \int_{R_n^+} \varphi(x) \prod_{i=1}^n j_{\gamma_i - \frac{1}{2}}(x_i y_i) d\mu_\gamma(x),$$

$$[F_{\mathrm{B}}^{-1}\ \varphi](x) = \int_{R_n^+} \varphi(y) \prod_{i=1}^n j_{\gamma_i - \frac{1}{2}}(x_i y_i) d\mu_\gamma(y)$$

where $c_\gamma = \prod_{i=1}^n \left[2^{\gamma_i - \frac{1}{2}} \Gamma(\gamma_i + \frac{1}{2})\right]^{-1}$ and $j_{\gamma_i - \frac{1}{2}}$ is the normalized Bessel function related to the Bessel function of the first kind by the formula $j_\gamma(r) = 2^\gamma \Gamma(\gamma + 1) J_\gamma(r) r^{-\gamma}$ [14]. However, the following equality for Fourier-Bessel transformation is true

$$F_{\mathrm{B}}(f * \varphi)(x) = F_{\mathrm{B}} f(x) F_{\mathrm{B}} \varphi(x).$$

## 2. The High Order Riesz Bessel Transforms $R_{\mathrm{B}_{x_i}}$ Associated with Bessel Generalized Shift Operator

In this section, we consider Bessel Generalized Shift Operator related to the Bessel operator. Then we give the Fourier-Bessel transformation of homogeneous polynomial which holds Bessel equations. Finally, we define the high order Riesz-Bessel transforms related to Bessel Generalized Shift Operator and so we show that high order Riesz-Bessel transforms hold the condition of classical Riesz transforms that is, these operators extend to high order Riesz-Bessel transforms [21].

It follows from the general theory of singular integrals that Riesz transforms are bounded on $L_{p,\gamma}(\mathbb{R}^n, d\mu_\gamma(x))$ for all $1 < p < \infty$. What is done in this paper is to extend this result to the context of Bessel theory where a similar operator is already defined. It has been noted that the difficulty arises in the application of the classical $L_p$-theory of Calderon-Zygmund, since Riesz transforms are singular integral operators. In this paper we describe how this theory can be adapted in Bessel setting and gives an $L_{p,\gamma}$-result for high order Riesz transforms for all $1 < p < \infty$.

$\Omega(y) = P_k(y)|y|^{-k}$, $K(y) = \Omega(y)|x|^{-n-2|\gamma|}$, the $P_k$ range over the homogeneous harmonic polynomials the latter arise in special case $k = 1$. Those for $k > 1$, we call the high order Riesz Bessel transform where we refer to $k$ as the degree of the high order Riesz Bessel transform [21]. They can also be characterized by their invariance properties.

Let $P_k$ be homogeneous polynomial of degree $k$ in $\mathbb{R}^n_+$. We shall say that $P$ is elliptic if $P(x)$ vanishes only at the origin. For any polynomial $P$ we consider also its corresponding differential polynomial. Thus if $P(x) = \sum a_\alpha x^\alpha$ we write $P(\frac{\partial}{\partial x}) = \sum a_\alpha (\frac{\partial}{\partial x})^\alpha$, where $(\frac{\partial}{\partial x})^\alpha = (\frac{\partial}{\partial x_1})^{\alpha_1} \dots (\frac{\partial}{\partial x_n})^{\alpha_n}$ and with the monomials $x^\alpha = x_1^{\alpha_1} \dots x_n^{\alpha_n}$ (see [20]).

**Theorem 2.1.** *Suppose that $P_k(x)$ is a homogeneous polynomial of degree $k$ and satisfies for Bessel operator $\mathrm{B}_{x_i}[P_k(x)] = 0$ then we have*

$$F_{\mathrm{B}}\left[P_k(x)e^{-|x|^2}\right](y) = 2^{-(|\gamma|+k+\frac{n}{2})} i^k P_k(y) e^{\frac{-|y|^2}{4}}.$$

**Lemma 2.2.** *Let $\theta = (\theta_1, \theta_2, \dots, \theta_n)$. Suppose that*

$$\int_{S^n_+} f(\theta) d\mu_\gamma(\theta) dS = 0$$

*and $\varphi$ is of Schwartz class $\mathcal{S}(\mathbb{R}^n_+)$ then we have the identity*

$$\lim\nolimits_{\varepsilon\to 0} \int_{\mathbb{R}^n} \frac{f(\frac{x}{|x|})}{|x|^{n+2|\gamma|-\varepsilon}} \varphi(x) d\mu_\gamma(x) = \lim\nolimits_{\varepsilon\to 0} \int_{|x|>\varepsilon} \frac{f(\frac{x}{|x|})}{|x|^{n+2|\gamma|}} \varphi(x) d\mu_\gamma(x).$$

*Proof.* The proof follows immediately from the representation.

$$\begin{aligned} \int_{\mathbb{R}^n_+} \frac{f(\frac{x}{|x|})}{|x|^{n+2|\gamma|-\varepsilon}} \varphi(x) d\mu_\gamma(x) &= \int_{|x|\le 1} \frac{f(\frac{x}{|x|})}{|x|^{n+2|\gamma|-\varepsilon}} [\varphi(x) - \varphi(0)] d\mu_\gamma(x) \\ &+ \int_{|x|>1} \frac{f(\frac{x}{|x|})}{|x|^{n+2|\gamma|-\varepsilon}} \varphi(x) d\mu_\gamma(x) \end{aligned}$$

□

The mean value theorem for Bessel differential operators is very convenient to obtain singular integral operators generated by a Bessel generalized shift operator. Therefore, we studied the mean value formula related to the Bessel generalized shift operator for the solutions of the boundary value problem for the Bessel operator $B_{x_j} u = 0$.

The Bessel generalized shift operator is one of the most important generalized shift operator on the half line $\mathbb{R}_+ = [0, +\infty)$, [8, 9, 14]. The Bessel generalized translation is used while studying various problems connected with Bessel operators . Fourier-Bessel harmonic analysis, i.e. the part of harmonic analysis addressing various problems on Bessel (Hankel) integral transforms, is closely connected with Bessel generalized shift operator.

## 3. The Mean Value Formula

In this section, we determine the mean value formula for Bessel generalized shift operator. Let $\mathbb{R}^n$ be $n$-dimensional Euclidian space and $x = (x_1, x_2, \ldots, x_n)$, $\xi = (\xi_1, \xi_2, \ldots, \xi_n)$ be vectors in $\mathbb{R}^n$. Then $x.\xi = x_1\xi_1 + \ldots + x_n\xi_n$. Denote, $\mathbb{R}^n_+ = \{x \in \mathbb{R}^n : x_i > 0, 1 \leq i \leq n\}$ and $S^n_+ = \{x \in \mathbb{R}^n_+ : |x| = 1\}$, $|\gamma| = \gamma_1 + \ldots + \gamma_n$. We assume that $D^+ \subset \mathbb{R}^n_+$ and $\Omega^+$ its boundary. In this paper, we are mainly concerned with the mean value theorem. Now, we relate this concept in the following theorem.

**Theorem 3.1.** *Let $S^{n-1}_+$ be unit sphere centered at origin, contained in $\mathbb{R}^n_+$ and $u$ be even regular solution with respect to $x_j, (j = 1, \ldots, n)$ of $\mathrm{B}_{x_j} u = 0$ . Then the following formula holds*

$$\int_{S^{n-1}_+} u(r\theta_1, \ldots, r\theta_n) \prod_{i=1}^n \theta_i^{2\gamma_i} dS = \frac{\prod_{i=1}^n \Gamma(\gamma_i + \frac{1}{2})}{2^{n-1}\Gamma(|\gamma| + \frac{n}{2})} u(0). \tag{3.1}$$

*Let $T^y$ be the Bessel generalized shift operator and $u$ be even regular solution of $\mathrm{B}_{x_j} u = 0$ at any interior point of the region $D^+ \subseteq \mathbb{R}^n_+$. Also, the following formula is valid*

$$\int_{S^{n-1}_+} T^y u(r\theta_1, r\theta_2, \ldots, r\theta_n) \prod_{i=1}^n \theta_i^{2\gamma_i} dS = \frac{\prod_{i=1}^n \Gamma(\gamma_i + \frac{1}{2})}{2^{n-1}\Gamma(|\gamma| + \frac{n}{2})} u(y). \tag{3.2}$$

*Proof.* We assume that the continuity of $u$ and $v$ in the closed region $\Omega^+ + D^+$, continuity of the first and second derivatives of $u$ and $v$ in $D^+$, together with continuity of the first derivatives of $u$ in $\Omega^+ + D^+$ and the second derivatives of $u$ in $D^+$. In addition, the existence of the integrals over $D^+$ are assumed in Green's formula. The most important tool for the potential theory is provided by this formula in $n$-dimensional bounded region $D^+$ with volume element $dg = dx_1 \ldots dx_n$ and its boundary $\Omega^+$, which we assume to be piecewise smooth, two functions $u$ and $v$ are related by the Green's formula

$$\int_{\Omega^+} \left( u\frac{\partial v}{\partial n} - v\frac{\partial u}{\partial n} \right) \prod_{i=1}^n x_i^{2\gamma_i} d\Omega \quad = \int\int_{D^+} [u\mathrm{B}_{x_j} v - v\mathrm{B}_{x_j} u] \prod_{i=1}^n x_i^{2\gamma_i} dg. \tag{3.3}$$

Let $P' = (0, \ldots, 0)$ be a boundary point of $\Omega^+$. We set

$$v = \left[ (n + 2|\gamma| - 2) r^{n+2|\gamma|-2} \right]^{-1} + w(r)$$

such that $r^2 = x_1^2 + x_2^2 + \ldots + x_n^2$ and $w \in \mathcal{C}^2(D^+)$ denote the set of two times continuously differentiable functions on $D^+$, also even function with respect to all $x_i$-variables. We apply formula (3.3) to the region $D^+ \setminus D_\varepsilon^+ \subset D^+$, where $D_\varepsilon^+$ is an upper half sphere centered at $P'$ of radius $\varepsilon$. Letting $\varepsilon \to 0$ and considering $\left[(n+2|\gamma|-2)r^{n+2|\gamma|-2}\right]^{-1}$ of the fundamental solution of $\mathrm{B}_{x_j}$, we deduce that

$$\begin{aligned}\int\int_{D^+\setminus D_\epsilon^+}[u\mathrm{B}_{x_j}v- \quad & v\mathrm{B}_{x_j}u]\prod_{i=1}^{n}x_i^{2\gamma_i}\,dx = \int_{\Omega}\left(u\tfrac{\partial v}{\partial n}-v\tfrac{\partial u}{\partial n}\right)\prod_{i=1}^{n}x_i^{2\gamma_i}\,d\Omega\\ & -\int_{\Omega\epsilon}\left(u\tfrac{\partial v}{\partial n}-v\tfrac{\partial u}{\partial n}\right)\prod_{i=1}^{n}x_i^{2\gamma_i}\,d\Omega.\end{aligned} \tag{3.4}$$

Hence, let us write the function $v$ in (3.4), then we have

$$\begin{aligned}\int\int_{D^+}[u\mathrm{B}_{x_j}v \quad & -v\mathrm{B}_{x_j}u]\prod_{i=1}^{n}x_i^{2\gamma_i}\,dx - \int\int_{D_\epsilon^+}[u\mathrm{B}_{x_j}w - w\mathrm{B}_{x_j}u]\prod_{i=1}^{n}x_i^{2\gamma_i}\,dx\\ & -\int\int_{D_\epsilon^+}\left[(n+2|\gamma|-2)r^{n+2|\gamma|-2}\right]^{-1}\mathrm{B}_{x_j}u\prod_{i=1}^{n}x_i^{2\gamma_i}\,dx\\ & = \int_{\Omega}\left(u\tfrac{\partial v}{\partial n}-v\tfrac{\partial u}{\partial n}\right)\prod_{i=1}^{n}x_i^{2\gamma_i}\,d\Omega + \int_{\Omega_\epsilon}u\frac{1}{r^{n+2|\gamma|-1}}\prod_{i=1}^{n}x_i^{2\gamma_i}\,d\Omega\\ & -\int_{\Omega_\epsilon}\left(u\tfrac{\partial w}{\partial n}-w\tfrac{\partial u}{\partial n}\right)\prod_{i=1}^{n}x_i^{2\gamma_i}\,d\Omega\\ & -\int_{\Omega_\epsilon}\left[(n+2|\gamma|-2)r^{n+2|\gamma|-2}\right]^{-1}\tfrac{\partial u}{\partial n}\prod_{i=1}^{n}x_i^{2\gamma_i}\,d\Omega\end{aligned} \tag{3.5}$$

and

$$\begin{aligned}\int\int_{D^+}[u\mathrm{B}_{x_j}v- \quad & v\mathrm{B}_{x_j}u]\prod_{i=1}^{n}x_i^{2\gamma_i}\,dg = \int_{\Omega}\left(u\tfrac{\partial v}{\partial n}-v\tfrac{\partial u}{\partial n}\right)\prod_{i=1}^{n}x_i^{2\gamma_i}\,d\Omega\\ & + \left[(n+2|\gamma|-2)r^{n+2|\gamma|-2}\right]^{-1}\int\int_{D_\epsilon^+}\mathrm{B}_{x_j}u\prod_{i=1}^{n}x_i^{2\gamma_i}\,dg\\ & - \left[(n+2|\gamma|-2)r^{n+2|\gamma|-2}\right]^{-1}\int_{\Omega_\epsilon}\frac{\partial u}{\partial n}\prod_{i=1}^{n}x_i^{2\gamma_i}\,d\Omega\\ & + \int_{\Omega_\epsilon}u\frac{1}{r^{n+2|\gamma|-1}}\prod_{i=1}^{n}x_i^{2\gamma_i}\,d\Omega.\end{aligned}$$

Therefore, we get

$$\begin{aligned}\int\int_{D^+}[u\mathrm{B}_{x_j}v - v\mathrm{B}_{x_j}u]\prod_{i=1}^{n}x_i^{2\gamma_i}\,dg \quad & = \int_{\Omega}\left(u\tfrac{\partial v}{\partial n}-v\tfrac{\partial u}{\partial n}\right)\prod_{i=1}^{n}x_i^{2\gamma_i}\,d\Omega\\ & + \int_{\Omega_\epsilon}u\frac{1}{r^{n+2|\gamma|-1}}\prod_{i=1}^{n}x_i^{2\gamma_i}\,d\Omega.\end{aligned}$$

Considering $\mathrm{B}_{x_j} v = \mathrm{B}_{x_j} w$, consequently, we have

$$\int\int_{D^+} [u\mathrm{B}_{x_j} v - v\mathrm{B}_{x_j} u] \prod_{i=1}^{n} x_i^{2\gamma_i} dg = \frac{\prod_{i=1}^{n} \Gamma(\gamma_i + \frac{1}{2})}{2^{n-1}\Gamma(|\gamma| + \frac{n}{2})} u(0) + \int_{\Omega} \left(u\frac{\partial v}{\partial n} - v\frac{\partial u}{\partial n}\right) \prod_{i=1}^{n} x_i^{2\gamma_i} d\Omega. \tag{3.6}$$

Here,

$$\int_{\Omega_\epsilon} \frac{1}{r^{n+2|\gamma|-1}} \prod_{i=1}^{n} x_i^{2\gamma_i} d\Omega = \frac{\prod_{i=1}^{n} \Gamma(\gamma_i + \frac{1}{2})}{2^{n-1}\Gamma(|\gamma| + \frac{n}{2})} = m(\Omega_\epsilon)$$

and

$$\lim_{\epsilon \to 0} \left[\int_{\Omega_\epsilon} u \frac{1}{r^{n+2|\gamma|-1}} \prod_{i=1}^{n} x_i^{2\gamma_i} d\Omega\right] = m(\Omega_\epsilon) u(0).$$

Let $K_R^+$ be a half sphere of radius $R$ centered at $P'$ contained in $D^+$ and $W_R^+$ be the surface of $K_R^+$. Suppose that

$$v = (n + 2|\gamma| - 2)^{-1} \left[r^{-(n+2|\gamma|-2)} - R^{-(n+2|\gamma|-2)}\right]. \tag{3.7}$$

Under our assumptions, we can rewrite the equality (3.6) as

$$\begin{aligned} \int_{K_R^+} \{u\mathrm{B}_{x_j} (n + 2|\gamma| - 2)^{-1} &\left[r^{-(n+2|\gamma|-2)} - R^{-(n+2|\gamma|-2)}\right] - v\mathrm{B}_{x_j} u\} \prod_{i=1}^{n} x_i^{2\gamma_i} dg \\ &= \int_{W_R^+} \left(u\frac{\partial v}{\partial n} - v\frac{\partial u}{\partial n}\right) \prod_{i=1}^{n} x_i^{2\gamma_i} dw \\ &+ \frac{\prod_{i=1}^{n} \Gamma(\gamma_i + \frac{1}{2})}{2^{n-1}\Gamma(|\gamma| + \frac{n}{2})} u(0). \end{aligned}$$

Since $R$ is constant and $\left[(n+2|\gamma|-2)r^{n+2|\gamma|-2}\right]^{-1}$ of the fundamental solution of $\mathrm{B}_{x_j}$, hence we have

$$\begin{aligned}
-\int_{K_R^+} v\mathrm{B}_{x_j} u \prod_{i=1}^n x_i^{2\gamma_i} dg &= \frac{\prod\limits_{i=1}^n \Gamma(\gamma_i+\frac{1}{2})}{2^{n-1}\Gamma(|\gamma|+\frac{n}{2})} u(0) \\
&+ \int_{W_R^+} \left\{ \left[u\tfrac{\partial (n+2|\gamma|-2)^{-1}}{\partial n} - (n+2|\gamma|-2)^{-1}\right] \right. \\
&\times \left. \left[r^{-(n+2|\gamma|-2)} - R^{-(n+2|\gamma|-2)}\right] \tfrac{\partial u}{\partial n} \right\} \prod_{i=1}^n x_i^{2\gamma_i} dw \\
&= \frac{\prod\limits_{i=1}^n \Gamma(\gamma_i+\frac{1}{2})}{2^{n-1}\Gamma(|\gamma|+\frac{n}{2})} u(0) \\
&+ \int_{W_R^+} u \tfrac{\partial (n+2|\gamma|-2)^{-1} r^{-(n+2|\gamma|-2)}}{\partial n} \prod_{i=1}^n x_i^{2\gamma_i} dw \\
&= \frac{\prod\limits_{i=1}^n \Gamma(\gamma_i+\frac{1}{2})}{2^{n-1}\Gamma(|\gamma|+\frac{n}{2})} u(0) - \frac{1}{R^{n+2|\gamma|-1}} \int_{W_R^+} u \prod_{i=1}^n x_i^{2\gamma_i} dw.
\end{aligned}$$

and

$$\begin{aligned}
\frac{2^{n-1}\Gamma(|\gamma|+\frac{n}{2})}{\prod\limits_{i=1}^n \Gamma(\gamma_i+\frac{1}{2}) R^{n+2|\gamma|-1}} \quad & \int_{W_R^+} u(x_1,\ldots,x_n) \prod_{i=1}^n x_i^{2\gamma_i} dw = u(0) + \\
& + \frac{2^{n-1}\Gamma(|\gamma|+\frac{n}{2})}{\prod\limits_{i=1}^n \Gamma(\gamma_i+\frac{1}{2})} \int_{K_R^+} v\mathrm{B}_{x_j} u \prod_{i=1}^n x_i^{2\gamma_i} dg.
\end{aligned} \tag{3.8}$$

Choose the function $v$ such that

$$v = c\frac{2^{n-1}\Gamma(|\gamma|+\frac{n}{2})}{\prod\limits_{i=1}^n \Gamma(\gamma_i+\frac{1}{2})} \left[(n+2|\gamma|-2)r^{n+2|\gamma|-2}\right]^{-1} + w(r), \tag{3.9}$$

where $w(r) \in \mathcal{C}^{(2)}(K_R^+)$, also even function with respect to all variables and $0 < r \leq R$. Then, under initial conditions

$$v|_{W_R^+} = \tfrac{\partial v}{\partial r}|_{W_R^+} = 0, \tag{3.10}$$

we obtain

$$\begin{aligned}
\int_{K_R^+} \left[u\mathrm{B}_{x_j} v - v\mathrm{B}_{x_j} u\right] \prod_{i=1}^n x_i^{2\gamma_i} dg &= \int_{W_R^+} \left(u\tfrac{\partial v}{\partial n} - v\tfrac{\partial u}{\partial n}\right) \prod_{i=1}^n x_i^{2\gamma_i} dw \\
&= \int_{W_R^+} u\tfrac{\partial v}{\partial n} \prod_{i=1}^n x_i^{2\gamma_i} dw \\
&= c\frac{2^{n-1}\Gamma(|\gamma|+\frac{n}{2})}{\prod\limits_{i=1}^n \Gamma(\gamma_i+\frac{1}{2})} \int_{W_R^+} u \frac{1}{r^{n+2|\gamma|-1}} \prod_{i=1}^n x_i^{2\gamma_i} dw \\
&= cu(0).
\end{aligned} \tag{3.11}$$

Since (3.11) is true for all $u$, we can take an arbitrary function that even with respect to all variables in $\mathcal{C}^{(2)}(K_R^+)$. Therefore, we assume that $u \in \mathcal{C}^{(2m+2)}(K_R^+)$ . For $\eta \leq m$, then we have the following identity

$$c\mathrm{B}_{x_j}^{\eta} u(0) = \int_{K_R^+} \left[\mathrm{B}_{x_j}^{\eta} u \mathrm{B}_{x_j} v - v\mathrm{B}_{x_j}^{\eta+1} u\right] \prod_{i=1}^{n} x_i^{2\gamma_i} dg \tag{3.12}$$

where $\mathrm{B}_{x_j}^{\eta}$ denotes the $\eta$-th order of $\mathrm{B}_{x_j}$, that is $\mathrm{B}_{x_j}^1 u = \mathrm{B}_{x_j} u$ , $\mathrm{B}_{x_j}^2 u = \mathrm{B}_{x_j}\mathrm{B}_{x_j} u$, etc. Now, let us define $\{v_\eta\}$ as a sequence of functions of the type (3.9). Then the following differential equation can be written

$$\mathrm{B}_{x_j} v_{\eta+1} = v_{\eta+1}^{''} + \tfrac{(n+2|\gamma|-1)}{r} v_{\eta+1}^{'} = v_\eta \quad (\eta = 0, 1, 2 \ldots, m). \tag{3.13}$$

By the initial conditions (3.10) and function

$$v_0 = \frac{2^{n-1}\Gamma(|\gamma| + \frac{n}{2})}{\prod\limits_{i=1}^{n} \Gamma(\gamma_i + \frac{1}{2})} (n + 2|\gamma| - 2)^{-1} \left[r^{-(n+2|\gamma|-2)} - R^{-(n+2|\gamma|-2)}\right], \tag{3.14}$$

it can be seen that the solution of equation (3.13) as follows

$$v_{\eta+1} = \left[(n + 2|\gamma| - 2) r^{n+2|\gamma|-2}\right]^{-1} \int_r^R \rho v_\eta(\rho) \left[\rho^{n+2|\gamma|-2} - r^{n+2|\gamma|-2}\right] d\rho.$$

Under the initial conditions, for each $\eta = 0, 1, 2 \ldots$, it can be easily found the solutions $v_0, v_1, \ldots$ corresponding to the constants $c_0, c_1, \ldots$ obtained as follows

$$c_\eta = \left(\tfrac{R}{2}\right)^{2\eta} \frac{\Gamma(|\gamma| + \frac{n}{2})}{\eta!\Gamma(\eta + |\gamma| + \frac{n}{2})} \quad (\eta = 0, 1, 2 \ldots).$$

Replacing $v$ by $v_\eta$ in (3.12), we see that

$$c_\eta \mathrm{B}_{x_j}^{\eta} u(0) = \int_{K_R^+} \left[\mathrm{B}_{x_j}^{\eta} u \mathrm{B}_{x_j} v_\eta - v_\eta \mathrm{B}_{x_j}^{\eta+1} u\right] \prod_{i=1}^{n} x_i^{2\gamma_i} dg.$$

Since $B_j v_\eta = v_{\eta-1}$, we obtain

$$c_\eta \mathrm{B}_{x_j}^{\eta} u(0) \quad = \quad \int_{K_R^+} \left[v_{\eta-1} \mathrm{B}_{x_j}^{\eta} u - v_\eta \mathrm{B}_{x_j}^{\eta+1} u\right] \prod_{i=1}^{n} x_i^{2\gamma_i} dg.$$

Summing all the equations corresponding to $\eta = 1, 2, \ldots, m$, then we conclude

$$\sum_{\eta=1}^{m} c_\eta \mathrm{B}_{x_j}^{\eta} u(0) = \int_{K_R^+} \left[v_0 \mathrm{B}_{x_j} u - v_m \mathrm{B}_{x_j}^{m+1} u\right] \prod_{i=1}^{n} x_i^{2\gamma_i} dg.$$

By (3.8) and the function $v_0$, then we get

$$\int_{K_R^+} v_0 \mathrm{B}_{x_j} u \prod_{i=1}^{n} x_i^{2\gamma_i} dg = \sum_{\eta=1}^{m} c_\eta \mathrm{B}_{x_j}^{\eta} u(0) + \int_{K_R^+} v_m \mathrm{B}_{x_j}^{m+1} u \prod_{i=1}^{n} x_i^{2\gamma_i} dg.$$

and

$$\frac{2^{n-1}\Gamma(|\gamma|+\frac{n}{2})}{\prod\limits_{i=1}^{n}\Gamma(\gamma_i+\frac{1}{2})}\int_{K_R^+} (n+2|\gamma|-2)^{-1}\left[r^{-(n+2|\gamma|-2)}-R^{-(n+2|\gamma|-2)}\right]$$
$$\times \mathrm{B}_{x_j}u\prod_{i=1}^{n}x_i^{2\gamma_i}dg=\sum_{\eta=1}^{m}\left(\frac{R}{2}\right)^{2\eta}\frac{\Gamma\big(|\gamma|+\frac{n}{2}\big)}{\eta!\Gamma\big(\eta+|\gamma|+\frac{n}{2}\big)}$$
$$\times \mathrm{B}_{x_j}^{\eta}u(0)+\int_{K_R^+}v_m\mathrm{B}_{x_j}^{m+1}u\prod_{i=1}^{n}x_i^{2\gamma_i}dg.$$

By the function $v$ in (3.7), we get

$$\frac{2^{n-1}\Gamma(|\gamma|+\dfrac{n}{2})}{\prod\limits_{i=1}^{n}\Gamma(\gamma_i+\frac{1}{2})}\int_{K_R^+}v\mathrm{B}_{x_j}u\prod_{i=1}^{n}x_i^{2\gamma_i}dg=\sum_{\eta=1}^{m}\left(\frac{R}{2}\right)^{2\eta}\frac{\Gamma(|\gamma|+\frac{n}{2})}{\eta!\Gamma(\eta+|\gamma|+\frac{n}{2})}\mathrm{B}_{x_j}^{\eta}u(0)$$
$$+\int_{K_R^+}v_m\mathrm{B}_{x_j}^{m+1}u\prod_{i=1}^{n}x_i^{2\gamma_i}dg \tag{3.15}$$

and from (3.8)

$$\left[c_{n,\gamma}R^{n+2|\gamma|-1}\right]^{-1}\int_{W_R^+}u(x)\prod_{i=1}^{n}x_i^{2\gamma_i}dw=u(0)+\sum_{\eta=1}^{m}\left(\frac{R}{2}\right)^{2\eta}\frac{\Gamma(|\gamma|+\frac{n}{2})}{\eta!\Gamma(\eta+|\gamma|+\frac{n}{2})}$$
$$\times\mathrm{B}_{x_j}^{\eta}u(0)+\int_{K_R^+}v_m\mathrm{B}_{x_j}^{m+1}u\prod_{i=1}^{n}x_i^{2\gamma_i}dg.$$

Thus we have

$$\left[c_{n,\gamma}R^{n+2|\gamma|-1}\right]^{-1}\int_{W_R^+}u(x)\prod_{i=1}^{n}x_i^{2\gamma_i}dw=\sum_{\eta=0}^{m}\left(\frac{R}{2}\right)^{2\eta}\frac{\Gamma(|\gamma|+\frac{n}{2})}{\eta!\Gamma(\eta+|\gamma|+\frac{n}{2})}$$
$$\times\mathrm{B}_{x_j}^{\eta}u(0)+\int_{K_R^+}v_m\mathrm{B}_{x_j}^{m+1}u\prod_{i=1}^{n}x_i^{2\gamma_i}dg, \tag{3.16}$$

where $c_0=1$, $\mathrm{B}_{x_j}^{0}u_0=u_0$ and $c_{n,\gamma}=\dfrac{\prod\limits_{i=1}^{n}\Gamma(\gamma_i+\frac{1}{2})}{2^{n-1}\Gamma(|\gamma|+\frac{n}{2})}$. This formula is valid for all arbitrary function $u\in K_R^+\subseteq D^+$ which is continuously differentiable function order $(2m+2)$ . Let us establish the equality (3.8) with respect to $T^yu(x)$, then we get

$$\frac{2^{n-1}\Gamma(|\gamma|+\frac{n}{2})}{\prod\limits_{i=1}^{n}\Gamma(\gamma_i+\frac{1}{2})R^{n+2|\gamma|-1}}\int_{W_R^+}T^yu(x)\prod_{i=1}^{n}x_i^{2\gamma_i}dw=u(y)$$
$$+\frac{2^{n-1}\Gamma(|\gamma|+\frac{n}{2})}{\prod\limits_{i=1}^{n}\Gamma(\gamma_i+\frac{1}{2})}\int_{K_R^+}v\mathrm{B}_{x_j}T^yu\prod_{i=1}^{n}x_i^{2\gamma_i}dg. \tag{3.17}$$

Setting $v = v_\eta$ and $u(x) = T^y u(x)$ in (3.12), then we obtain

$$c\mathrm{B}^{\eta}_{x_j} u(y) = \int_{K_R^+} \left[\mathrm{B}^{\eta}_{x_j} T^y u(x) \mathrm{B}_{x_j} v - v \mathrm{B}^{\eta+1}_{x_j} T^y u(x)\right] \prod_{i=1}^{n} x_i^{2\gamma_i} dg$$

$$c_\eta \mathrm{B}^{\eta}_{x_j} u(y) = \int_{K_R^+} \left[\mathrm{B}^{\eta}_{x_j} T^y u(x) \mathrm{B}_{x_j} v_\eta - v_\eta \mathrm{B}^{\eta+1}_{x_j} T^y u(x)\right] \prod_{i=1}^{n} x_i^{2\gamma_i} dg.$$

If we sum all the equations corresponding to $\eta = (1, 2, \dots, m)$, we conclude

$$\sum_{\eta=1}^{m} c_\eta \mathrm{B}^{\eta}_{x_j} u(y) = \int_{K_R^+} \left[v_0 \mathrm{B}_{x_j} T^y u - v_m \mathrm{B}^{m+1}_{x_j} T^y u\right] \prod_{i=1}^{n} x_i^{2\gamma_i} dg.$$

By (3.8) and the function $v_0$, then we get

$$\int_{K_R^+} v_0 \mathrm{B}_{x_j} T^y u(x) \prod_{i=1}^{n} x_i^{2\gamma_i} dg = \sum_{\eta=1}^{m} \left(\tfrac{R}{2}\right)^{2\eta} \frac{\Gamma(|\gamma| + \frac{n}{2})}{\eta! \Gamma(\eta + |\gamma| + \frac{n}{2})} \mathrm{B}^{\eta}_{x_j} u(y) + \int_{W_R^+} v_m \mathrm{B}^{m}_{x_j} T^y u \prod_{i=1}^{n} x_i^{2\gamma_i} dw \tag{3.18}$$

and

$$\frac{2^{n-1}\Gamma(|\gamma| + \frac{n}{2})}{\prod_{i=1}^{n} \Gamma(\gamma_i + \frac{1}{2}) R^{n+2|\gamma|-1}} \int_{K_R^+} v \mathrm{B}_{x_j} T^y u(x) \prod_{i=1}^{n} x_i^{2\gamma_i} dg = \sum_{\eta=1}^{m} \left(\tfrac{R}{2}\right)^{2\eta} \frac{\Gamma(|\gamma| + \frac{n}{2})}{\eta! \Gamma(\eta + |\gamma| + \frac{n}{2})} \times \mathrm{B}^{\eta}_{x_j} u(y) + \int_{K_R^+} v_m \mathrm{B}^{m+1}_{x_j} T^y u(x) \prod_{i=1}^{n} x_i^{2\gamma_i} dg \tag{3.19}$$

$$\left[c_{n,\gamma} R^{n+2|\gamma|-1}\right]^{-1} \int_{W_R^+} T^y u(x) \prod_{i=1}^{n} x_i^{2\gamma_i} dw_R = u(y) + \sum_{\eta=1}^{m} \left(\tfrac{R}{2}\right)^{2\eta} \frac{\Gamma(|\gamma| + \frac{n}{2})}{\eta! \Gamma(\eta + |\gamma| + \frac{n}{2})} \times \mathrm{B}^{\eta}_{x_j} u(y) + \int_{K_R^+} v_m \mathrm{B}^{m+1}_{x_j} T^y u(x) \prod_{i=1}^{n} x_i^{2\gamma_i} dg \tag{3.20}$$

and

$$\left[c_{n,\gamma} R^{n+2|\gamma|-1}\right]^{-1} \int_{W_R^+} T^y u(x) \prod_{i=1}^{n} x_i^{2\gamma_i} dw_R = \sum_{\eta=0}^{m} \left(\tfrac{R}{2}\right)^{2\eta} \frac{\Gamma(|\gamma| + \frac{n}{2})}{\eta! \Gamma(\eta + |\gamma| + \frac{n}{2})} \mathrm{B}^{\eta}_{x_j} u(y) + \int_{K_R^+} v_m \mathrm{B}^{m+1}_{x_j} T^y u(x) \prod_{i=1}^{n} x_i^{2\gamma_i} dg \tag{3.21}$$

where $\mathrm{B}^{0}_{x_j} u = u$ and $c_0 = 1$ . We note that the second integral on the right-hand side in (3.21) tends to zero for $m \to \infty$ then we have the result

$$\frac{2^{n-1}\Gamma(|\gamma| + \frac{n}{2})}{\prod_{i=1}^{n} \Gamma(\gamma_i + \frac{1}{2}) R^{n+2|\gamma|-1}} \int_{S_{n-1}^+} T^y u(R\theta_1, \dots, R\theta_n) R^{2|\gamma|} \prod_{i=1}^{n} \theta_i^{2\gamma_i} R^{n-1} dS = \sum_{\eta=0}^{\infty} \left(\tfrac{R}{2}\right)^{2\eta} \frac{\Gamma(|\gamma| + \frac{n}{2})}{\eta! \Gamma(\eta + |\gamma| + \frac{n}{2})} \mathrm{B}^{\eta}_{x_j} u(y)$$

where $x = R\theta$ and $dw_R = R^{n-1}dS$ . Hence

$$\int_{S_+^{n-1}} T^y u(R\theta_1, R\theta_2, \dots, R\theta_n) \prod_{i=1}^n \theta_i^{2\gamma_i} dS = \frac{\prod_{i=1}^n \Gamma(\gamma_i + \frac{1}{2})}{2^{n-1}\Gamma(|\gamma| + \frac{n}{2})} u(y).$$

Thus the proof is completed. □

*Proof of Theorem 2.1.* We shall need the Fourier-Bessel transforms of the function $e^{-\alpha|x|^2}$. The identity to be proved can be rewritten as

$$F_{\mathrm{B}}\big(e^{-\alpha|x|^2}\big)(y) \quad = c_\gamma \int_{\mathbb{R}_+^n} e^{-\alpha|x|^2} \prod_{i=1}^n j_{\gamma_i - \frac{1}{2}}(x_i y_i) d\mu_\gamma(x). \tag{3.22}$$

Let $\gamma > -1, \alpha > 0$ and $J_\gamma(br)$ be Bessel function. We recall that

$$\int_0^\infty e^{-\alpha r^2} r^{\gamma+1} J_\gamma(br) dr = b^\gamma (2\alpha)^{-\gamma-1} e^{\frac{-b^2}{4\alpha}},$$

then we get

$$J_{\gamma_1 - \frac{1}{2}}(x_1 y_1) = \big[2^{\gamma_1 - \frac{1}{2}} \Gamma\big(\gamma_1 + \frac{1}{2}\big)\big]^{-1} (x_1 y_1)^{\gamma_1 - \frac{1}{2}} j_{\gamma_1 - \frac{1}{2}}(x_1 y_1)$$

$$\vdots$$

$$J_{\gamma_n - \frac{1}{2}}(x_n y_n) = \big[2^{\gamma_n - \frac{1}{2}} \Gamma\big(\gamma_n + \frac{1}{2}\big)\big]^{-1} (x_n y_n)^{\gamma_n - \frac{1}{2}} j_{\gamma_n - \frac{1}{2}}(x_n y_n)$$

where $J_\gamma(r) = [2^\gamma \Gamma(\gamma + 1)]^{-1} r^\gamma j_\gamma(r)$. We may now calculate each integrals in equation (3.22).

$$\begin{aligned} I_1 \quad &= \int_0^\infty e^{-\alpha x_1{}^2} j_{\gamma_1 - \frac{1}{2}}(x_1 y_1) x_1^{2\gamma_1} dx_1 \\ &= \int_0^\infty e^{-\alpha r^2} J_{\gamma_1 - \frac{1}{2}}(rs) r^{2\gamma_1} r^{-\gamma_1 + \frac{1}{2}} s^{-\gamma_1 + \frac{1}{2}} dr \\ &= 2^{\gamma_1 - \frac{1}{2}} \Gamma\big(\gamma_1 + \tfrac{1}{2}\big) (2\alpha)^{-\gamma_1 - \frac{1}{2}} e^{-\frac{s^2}{4\alpha}} \end{aligned} \tag{3.23}$$

where $x_1 = r$ and $y_1 = s$. By the similar way, we obtain

$$I_2 \quad = \int_0^\infty e^{-\alpha x_2{}^2} j_{\gamma_2 - \frac{1}{2}}(x_2 y_2) x_2^{2\gamma_2} dx_2 = \frac{2^{\gamma_2 - \frac{1}{2}} \Gamma\big(\gamma_2 + \frac{1}{2}\big)}{(2\alpha)^{\gamma_2 + \frac{1}{2}}} e^{-\frac{y_2^2}{4\alpha}}, \tag{3.24}$$

$$I_n \quad = \int_0^\infty e^{-\alpha x_n{}^2} j_{\gamma_n - \frac{1}{2}}(x_n y_n) x_2^{2\gamma_n} dx_n = \frac{2^{\gamma_n - \frac{1}{2}} \Gamma\big(\gamma_n + \frac{1}{2}\big)}{(2\alpha)^{\gamma_n + \frac{1}{2}}} e^{-\frac{y_n^2}{4\alpha}}. \tag{3.25}$$

If we write these equalities (3.23)-(3.25) in (3.22) then we have

$$F_{\mathrm{B}}\big(e^{-\alpha|x|^2}\big)(y) = e^{-\frac{|y|^2}{4\alpha}} (2\alpha)^{\frac{-2|\gamma|-n}{2}}.$$

Letting $\alpha \to 1$ and $y \to -2y$, then we get

$$F_{\mathrm{B}}\big(e^{-\alpha|x|^2}\big)(y) = e^{-|y|^2} 2^{-n}.$$

So, we can easily obtain the identity

$$F_{\mathrm{B}}\big(e^{-|x|^2}\big)(y) \quad = c_\gamma \int_{\mathbb{R}^n_+} \prod_{i=1}^{n} j_{\gamma_i-\frac{1}{2}}(x_i 2y_i) e^{-|x|^2} d\mu_\gamma(x).$$

By the similar arguments, we can see that

$$I_1^{'} = \int_0^\infty e^{-{x_1}^2} j_{\gamma_1-\frac{1}{2}}(x_1 2y_1) x_1^{2\gamma_1} dx_1,$$

$$I_1^{'} = \int_0^\infty e^{-{r^{'}}^2} j_{\gamma_1-\frac{1}{2}}(2r^{'}s^{'}) {r^{'}}^{2\gamma_1} dr^{'}.$$

Hence, we get

$$I_1^{'} = 2^{-1}\Gamma\big(\gamma_1 + \tfrac{1}{2}\big) e^{-{s_1^{'}}^2}$$

and

$$I_n^{'} \quad = \int_0^\infty e^{-{r^{'}}^2} j_{\gamma_n-\frac{1}{2}}(2r^{'}s^{'}) {r^{'}}^{2\gamma_n} dr^{'} = 2^{-1} e^{-{s_n^{'}}^2} \Gamma\big(\gamma_n + \tfrac{1}{2}\big).$$

Consequently, we deduce

$$F_{\mathrm{B}}\big(e^{-\alpha|x|^2}\big)(y) = 2^{-n} \prod_{i=1}^{n} \Gamma\big(\gamma_i + \tfrac{1}{2}\big) e^{-|y|^2}. \tag{3.26}$$

When we apply the Bessel differential operator $P_k(\mathrm{B}_{t_1}, \mathrm{B}_{t_2}, \ldots, \mathrm{B}_{t_{n-1}}, \mathrm{B}_{t_n})$ to both sides of the identity (3.26), then we see

$$\begin{aligned} F_{\mathrm{B}}\big[P_k(x)e^{-|x|^2}\big](t) \quad &= \int_0^\infty P_k(x) e^{-{x_1}^2} j_{\gamma_1-\frac{1}{2}}(2x_1t_1) x_1^{2\gamma_1} \\ &\times \ldots \int_0^\infty e^{-{x_n}^2} j_{\gamma_n-\frac{1}{2}}(2x_nt_n) x_n^{2\gamma_n} dx \\ &= Q(t) 2^{-n} \prod_{i=1}^{n} \Gamma\big(\gamma_i + \tfrac{1}{2}\big) e^{-|t|^2} \end{aligned}$$

where $Q(t)$ is a polynomial. The problem is therefore to show that $Q(t) = P_k(iy)$. Now using the identity

$$j_{\gamma-\frac{1}{2}}(r) = \frac{\Gamma\big(\gamma + \frac{1}{2}\big)}{\Gamma(\gamma)\Gamma(\frac{1}{2})} \int_0^\pi e^{ir\cos\alpha} (\sin\alpha)^{2\gamma-1} d\alpha$$

then we have

$$\begin{aligned}
Q(t) &= c_1\int_{\mathbb{R}^n_+} P_k(x)e^{|t|^2-|x|^2}\prod_{i=1}^n j_{\gamma_i-\frac{1}{2}}(2x_it_i)d\mu_\gamma(x)\\
&= c_1\int_{\mathbb{R}^n_+} P_k(x)e^{t_1^2+\ldots+t_n^2-(x_1^2+\ldots+x_n^2)}\prod_{i=1}^n j_{\gamma_i-\frac{1}{2}}(2x_it_i)d\mu_\gamma(x)\\
&= c_2\int_{\mathbb{R}^n_+} P_k(x)e^{|t|^2-|x|^2}\Big(\int_0^\pi e^{2ix_1t_1\cos w}\sin^{2\gamma_1-1}wdw\\
&\quad \ldots\int_0^\pi e^{2ix_nt_n\cos w}\sin^{2\gamma_n-1}wdw\Big)d\mu_\gamma(x)\\
&= c_2\int_{\mathbb{R}^n_+} P_k(x)d\mu_\gamma(x)\int_0^\pi\ldots\int_0^\pi\Big(e^{t_1^2-x_1^2-2ix_1t_1\cos w}\sin^{2\gamma_1-1}w\Big)\\
&\quad \ldots\Big(e^{t_n^2-x_n^2-2ix_nt_n\cos w}\sin^{2\gamma_n-1}w\Big)dw = 2^n\int_{\mathbb{R}^n_+}P_k(x)[T^{-it}(e^{-|x|^2})]d\mu_\gamma(x)
\end{aligned}$$

where $c_1 = 2^n\big[\prod_{i=1}^n\Gamma\big(\gamma_i+\frac{1}{2}\big)\big]^{-1}$ and $c_2 = \prod_{i=1}^n 2^n\Gamma\big(\gamma_i+\frac{1}{2}\big)\big(\pi^{\frac{n}{2}}\Gamma(\gamma_i)\big)^{-1}$. Replacing $t$ by $-it$ and the properties of $T^y$, we obtain

$$Q(-it) = c_1\int_{\mathbb{R}^n_+} T^{-t}[P_k(x)e^{-|x|^2}]d\mu_\gamma(x). \tag{3.27}$$

Taking the change of variables $x\to r\theta$ for $0<r<\infty$ and $\theta\in\mathcal{S}^{n-1}_+$ and applying the polar coordinates in (3.27), this gives the identity

$$Q(-it) = c_1\int_0^\infty r^{2|\gamma|+n-1}\Big(\int_{\mathcal{S}^{n-1}_+} T^{-t}[P_k(r\theta)]\prod_{i=1}^n\theta_i^{2\gamma_i}\Big)e^{-r^2}drd\theta. \tag{3.28}$$

We must calculate the integral $\int_0^\infty r^{2|\gamma|+n-1}e^{-r^2}dr$. Letting $r^2=u$, then we get

$$\int_0^\infty r^{2|\gamma|+n-1}e^{-r^2}dr = \tfrac{1}{2}\int_0^\infty u^{|\gamma|+\frac{n}{2}-1}e^{-u}dr = \tfrac{1}{2}\Gamma\big(|\gamma|+\tfrac{n}{2}\big).$$

Applying Theorem 3.1 to (3.28), we deduce $Q(-it)=P_k(t)$ and

$$\begin{aligned}
F_{\mathrm{B}}\big[P_k(x)e^{-\alpha|x|^2}\big](t) &= c_\gamma\int_{\mathbb{R}^n_+}P_k(x)e^{-|x|^2}\prod_{i=1}^n j_{\gamma_i-\frac{1}{2}}(2x_it_i)d\mu_\gamma(x)\\
&= 2^{-(|\gamma|+\frac{n}{2})}P_k(it)e^{-|t|^2}.
\end{aligned}$$

Since $P_k$ is homogeneous, we have the identity

$$F_{\mathrm{B}}\big[P_k(x)e^{-|x|^2}\big](y) = 2^{-(|\gamma|+k+\frac{n}{2})}i^kP_k(y)e^{\frac{-|y|^2}{4}}$$

and so we obtain the desired conclusion. □

We come now to what has been our main goal in this paper.

**Theorem 3.2.** *Let $P_k$ be homogeneous harmonic polynomial of degree $k$. Then we get*

$$F_{\mathrm{B}}\big[\ p.v\ \frac{P_k(x)}{|x|^{k+n+2|\gamma|}}\big](y) = 2^{\frac{-n-2|\gamma|}{2}} i^k \frac{\Gamma(\frac{k}{2})}{\Gamma\big(\frac{k+n+2|\gamma|}{2}\big)} \frac{P_k(y)}{|y|^k}.$$

*Proof.* Consider the identity

$$F_{\mathrm{B}}[f(\alpha x)](t) = \alpha^{-n-2|\gamma|} F_{\mathrm{B}}[f(x)](\tfrac{t}{\alpha}). \tag{3.29}$$

By (3.29) and Theorem 2.1, we have

$$F_{\mathrm{B}}[P_k(x)e^{-\alpha|x|^2}](y) = (2\alpha)^{-(k+|\gamma|+\frac{n}{2})} i^k e^{\frac{-|y|^2}{4\alpha}} P_k(y).$$

Assume in addition that $\varphi \in \mathcal{S}(\mathbb{R}^n_+)$. Then

$$\int_{\mathbb{R}^n_+} P_k(x)e^{-\alpha|x|^2}\varphi(x)d\mu_\gamma(x) \quad = (2\alpha)^{-(k+|\gamma|+\frac{n}{2})} i^k \int_{\mathbb{R}^n_+} P_k(x)e^{-\frac{|x|^2}{4\alpha}}\varphi(x)d\mu_\gamma(x).$$

We now integrate both sides of the above with respect to $\alpha$ having multiplied the equation by a suitable power of $\alpha$ ($\alpha^{\beta-1},\ 2\beta = k+n+2|\gamma|-\varepsilon$). We obtain

$$\int_0^\infty \big[\int_{\mathbb{R}^n_+} \quad P_k(x)e^{-\alpha|x|^2}\alpha^{\beta-1}\varphi(x)d\mu_\gamma(x)\big]d\alpha = \int_0^\infty \big[(2\alpha)^{-(k+|\gamma|+\frac{n}{2})} i^k \alpha^{\beta-1}\big]$$
$$\times\big[\int_{\mathbb{R}^n_+} P_k(x)e^{-\frac{|x|^2}{4\alpha}}\varphi(x)d\mu_\gamma(x)\big]d\alpha.$$

If we use the fact that

$$\int_0^\infty e^{-u}\big(\tfrac{u}{|x|^2}\big)^{\beta-1}\tfrac{du}{|x|^2} = \Gamma\big(\beta\big)|x|^{-2\beta}$$

we get

$$\int_0^\infty \big[\int_{\mathbb{R}^n_+} P_k(x)e^{-\alpha|x|^2}\alpha^{\beta-1}\varphi(x)d\mu_\gamma(x)\big]d\alpha \quad = \Gamma\big(\beta\big)\int_{\mathbb{R}^n_+} \frac{P_k(x)}{|x|^{2\beta}}\varphi(x)d\mu_\gamma(x). \tag{3.30}$$

The corresponding integration for the right side gives

$$\int_0^\infty \quad \big\{(2\alpha)^{-(k+|\gamma|+\frac{n}{2})} i^k \alpha^{\beta-1}\big[\int_{\mathbb{R}^n_+} P_k(x)e^{-\frac{|x|^2}{4\alpha}}\varphi(x)d\mu_\gamma(x)\big]\big\}d\alpha$$
$$= \int_0^\infty 2^{-(k+|\gamma|+\frac{n}{2})} i^k \alpha^{-(\frac{k+\varepsilon}{2}+1)} \int_{\mathbb{R}^n_+} \big[P_k(x)e^{-\frac{|x|^2}{4\alpha}}\varphi(x)d\mu_\gamma(x)\big]d\alpha$$
$$= 2^{-(k+|\gamma|+\frac{n}{2})} \int_0^\infty \int_{\mathbb{R}^n_+} \alpha^{(\frac{-k-\varepsilon}{2}-1)} P_k(x)e^{-\frac{|x|^2}{4\alpha}} i^k \varphi(x)d\mu_\gamma(x)d\alpha.$$

Letting $\alpha \to \frac{1}{4\alpha}$, we obtain

$$\int_0^\infty \Big\{(2\alpha)^{-(k+|\gamma|+\frac{n}{2})} i^k \alpha^{\beta-1}\Big[\int_{\mathbb{R}^n_+} P_k(x) e^{-\frac{|x|^2}{4\alpha}} \varphi(x) d\mu_\gamma(x)\Big]\Big\} d\alpha$$
$$= 2^{-(|\gamma|+\frac{n}{2}-\varepsilon)} \Gamma\big(\tfrac{k+\varepsilon}{2}\big) i^k \int_{\mathbb{R}^n_+} \frac{P_k(x)}{|x|^{k+\varepsilon}} \varphi(x) d\mu_\gamma(x). \tag{3.31}$$

Thus (3.30) and (3.31) lead to the identity

$$\int_{\mathbb{R}^n_+} \frac{P_k(x)}{|x|^{2\beta}} \varphi(x) d\mu_\gamma(x) = 2^{-\frac{2|\gamma|+n}{2}+\varepsilon} i^k \frac{\Gamma\big(\frac{k+\varepsilon}{2}\big)}{\Gamma(\beta)} \int_{\mathbb{R}^n_+} \frac{P_k(x)}{|x|^{k+\varepsilon}} \varphi(x) d\mu_\gamma(x).$$

By Lemma 2.2, we have therefore concluded the proof of the theorem

$$F_B\Big[\, p.v\ \frac{P_k(x)}{|x|^{k+n+2|\gamma|-\varepsilon}}\Big](y) = 2^{-\frac{2|\gamma|+n}{2}} i^k \frac{\Gamma\big(\frac{k}{2}\big)}{\Gamma\big(\frac{k+n+2|\gamma|}{2}\big)} \frac{P_k(y)}{|y|^k}$$

□

**Definition 3.3.** Let $T^y$ be the Bessel Generalized Shift Operator and let $f$ be a Schwartz function on $\mathbb{R}^n_+$. We define the high order Riesz Bessel transforms $R^{(k)}_{\mathrm{B}_{x_i}}$ of order $k$ with respect to Bessel Generalized Shift Operator as

$$R^{(k)}_{\mathrm{B}}(f)(x) = c_k(n,\gamma)\Big[\, p.v\ \frac{P_k(y)}{|y|^{n+k+2|\gamma|}} * f\Big](x)$$
$$= c_k(n,\gamma) \lim_{\varepsilon\to 0} \int_{0<\varepsilon<|x|} \frac{P_k(y)}{|y|^{k+n+2|\gamma|}} T^y f(x) d\mu_\gamma(y),$$

where $c_k(n,\gamma) = 2^{\frac{n+2|\gamma|}{2}} \Gamma\big(\frac{n+k+2|\gamma|}{2}\big)\big[\Gamma\big(\frac{k}{2}\big)\big]^{-1} (k = 1, 2, \ldots, n)$ and $P_k(x)$ is a homogeneous polynomial of degree $k$ in $\mathbb{R}^n_+$ which satisfies $\mathrm{B}_{x_i} P_k = 0$.

By Theorem 3.2, we conclude

$$F_{\mathrm{B}}[R^{(k)}_{\mathrm{B}}(f)](\xi) = i^k P_k(\xi)|\xi|^{-k} F_{\mathrm{B}}[f](\xi). \tag{3.32}$$

One of the important applications of the high order Riesz transforms is that they can be used to mediate between various combinations of partial derivatives of a function. We shall here content ourselves with two very simple illustrations, which examples have an interest on their own and have already the characteristic features of a general type of estimate which can be made in the theory of elliptic differential operators.

**Proposition 3.4.** *Suppose $f$ is a class of $\mathcal{S}(\mathbb{R}^n_+)$ and has compact support. Let $\mathrm{B}_{x_i} f$ be the Bessel differential operator. Then we have the a priori bound*

$$||\partial_{x_i}\partial_{x_k} f||_{p,\gamma} \le A_p ||\mathrm{B}_{x_i} f||_{p,\gamma}$$

*with $A_p$ independent of $f$.*

In (3.32), if we take $k = 1$. Then this proposition is an immediate consequence of the $L_{p,\gamma}$ boundedness of the Riesz Bessel transforms generated by Bessel Generalized Shift Operator and the identity

$$\partial_{x_i}(\partial_{x_k} f) = -R_{\mathrm{B}_{x_i}} R_{B_{x_k}} \mathrm{B}_{x_i} f. \tag{3.33}$$

To prove (3.33) we use the Fourier-Bessel transform. Thus if $F_B[f](x)$ is the Fourier-Bessel transform of $f$, then the Fourier-Bessel transform of $\partial_{x_i} f$ is

$$F_{\mathrm{B}}[\partial_{x_i} f](y) = -x_i F_{\mathrm{B}}[f](y),$$

and

$$F_{\mathrm{B}}[\mathrm{B}_{x_i} f](y) = -|x|^2 F_{\mathrm{B}}[f](y)$$

and so

$$\begin{aligned} F_{\mathrm{B}}\left[\partial_{x_i}(\partial_{x_k} f)\right](y) &= -\big(\frac{x_i}{|x|}\big)\big(\frac{x_k}{|x|}\big)|x|^2 F_{\mathrm{B}}[f](y) \\ &= -F_{\mathrm{B}}[R_{\mathrm{B}_{x_i}} R_{B_{x_k}} \mathrm{B}_{x_i} f](y) \end{aligned}$$

which gives the (3.33). Thus we have

$$||\partial_{x_i}(\partial_{x_k} f)||_p = ||R_{\mathrm{B}_{x_i}} R_{\mathrm{B}_{x_k}} \mathrm{B}_{x_i} f||_p \leq C||\mathrm{B}_{x_i} f||_p.$$

By using the Fourier Bessel transformations, we have

$$P(|y|^2) F_{\mathrm{B}}[(\mathrm{B}_{x_i}) f](y) \quad = -|y|^2 F_{\mathrm{B}}[P(\mathrm{B}_{x_i}) f](y)$$

(see [5]).

**Corollary 3.5.** *Suppose $P_k$ is a homogeneous elliptic polynomial of degree $k$ and $f$ is $k$-times continuously differentiable with compact support. Then we have the priori estimate*

$$\big|\big|\mathrm{B}_{x_i} f\big|\big|_{p,\gamma} \leq A_p \big|\big|P_k(\mathrm{B}_{x_i}) f\big|\big|_{p,\gamma}, \qquad 1 < p < \infty.$$

To prove this inequality, we note that the following relation between Fourier Bessel transform of $\mathrm{B}_{x_i} f$ and $P_k(\mathrm{B}_{x_i}) f$ holds

$$P_k(|y|^2) F_{\mathrm{B}}\left[\mathrm{B}_{x_i} f\right](y) = -|y|^2 F_{\mathrm{B}}\left[P_k(\mathrm{B}_{x_i}) f\right](y).$$

Since $P_k(|y|^2)$ is non-vanishing except the origin and let $-\frac{|y|^2}{P_k(|y|^2)}$ be homogeneous of degree zero and indefinitely differentiable on the unit sphere. Then we get

$$\mathrm{B}_{x_i} f = R_{\mathrm{B}_{x_i}}\big(P_k(\mathrm{B}_{x_i}) f\big).$$

We also have the following $L_{p,\gamma}$ boundedness of the high order Riesz Bessel transform.

**Theorem 3.6.** *The high order Riesz-Bessel transforms generated by Bessel Generalized Shift Operator are bounded operator from $L_{p,\gamma}(\mathbb{R}^n_+)$ into itself for all $1 < p < \infty$*

$$||R_{\mathrm{B}_{x_i}} f||_{p,\gamma} \leq A_p ||f||_{p,\gamma}.$$

I. Ekincioglu
Department of Mathematics
Dumlupınar University
Kütahya, Turkey
e-mail: ismail.ekincioglu@dpu.edu.tr

H.H. Sayan
Technology Faculty
Department of Electric-Electronic Engineering
Gazi University
Ankara, Turkey
e-mail: hsayan@gazi.edu.tr

C. Keskin
Department of Mathematics
Dumlupınar University
Kütahya, Turkey
e-mail: cansu.keskin@dpu.edu.tr